\newtheorem{thm}{Theorem}

\newtheorem{prop}{Proposition}

\newtheorem{lemma}{Lemma}

\newtheorem{ex}{Example}

\PassOptionsToPackage{usenames,dvipsnames}{xcolor}

\documentclass[11pt]{article}
\usepackage{amsfonts}
\usepackage{mathtools}
\mathtoolsset{showonlyrefs}

\usepackage{blindtext}
\usepackage{tikz}

\usepackage{amsmath}
\usepackage{graphicx}

\newdimen\slantmathcorr
\def\oversl#1{
\setbox0=\hbox{$#1$}
\slantmathcorr=\wd0
\hskip 0.2\slantmathcorr \overline{\hbox to 0.8\wd0{%
\vphantom{\hbox{$#1$}}}}
\hskip-\wd0\hbox{$#1$}
}

\def\undersl#1{
\setbox0=\hbox{$#1$}
\slantmathcorr=\wd0
\underline{\hbox to 0.8\wd0{%
\vphantom{\hbox{$#1$}}}}
\hskip-0.8\wd0\hbox{$#1$}
}

\newcommand{\be}{\begin{equation}}
\newcommand{\ee}{\end{equation}}
\newcommand{\ben}{\begin{equation*}}
\newcommand{\een}{\end{equation*}}
\newcommand{\ba}{\begin{aligned}}
\newcommand{\ea}{\end{aligned}}

\newcommand{\geps}{\epsilon}
\newcommand{\gep}{w}

\newcommand{\rmi}{{\rm i}}

\newcommand{\wh}{\widehat}

\newcommand{\whH}{\widehat{H}}
\newcommand{\whV}{\widehat{V}}
\newcommand{\whL}{\widehat{L}}
\newcommand{\whX}{\widehat{X}}

\newcommand{\whK}{\widehat{K}}

\newcommand{\whk}{\widehat{\kappa}}

\newcommand{\whe}{\widehat{e}}

\newcommand{\Xbar}{\overline{X}}

\newcommand{\RbarG}{\overline{R}^\Gamma}

\newcommand{\wbar}{\overline{w}}

\newcommand{\nbar}{\overline{n}}

\newcommand{\R}{\mathbb{R}}

\newcommand{\tG}{\tilde G}
\newcommand{\Pu}{P^{(u)}}
\newcommand{\PuG}{P^{(u,\Gamma)}}

\newcommand{\EuG}{E^{(u,\Gamma)}}

\newcommand{\ga}{\alpha}
\newcommand{\gl}{\lambda}
\newcommand{\gk}{\kappa}
\newcommand{\gz}{\zeta}
\newcommand{\gb}{\beta}

\newcommand{\gd}{\delta}

\newcommand{\tK}{\tilde K}

\newcommand{\RG}{{R}^\Gamma}
\newcommand{\RGu}{{R}^{\Gamma,u}}
\newcommand{\mcRG}{\overline{\mathcal R}^\Gamma}

\newcommand{\mcT}{\mathcal T}

\newcommand{\Xubar}{\undersl{X}}
\newcommand{\whn}{\widehat{n}}

\newcommand{\gepsbar}{\overline{\geps}}
\newcommand{\ebar}{\overline{e}}

\newcommand{\tu}{\tau_u}

\newcommand{\tuG}{\tau^\Gamma_u}
\newcommand{\suG}{\sigma^\Gamma_u}
\newcommand{\tv}{\tau_v}
\newcommand{\too}{\tau_0}

\newcommand{\tz}{\tau_z}
\newcommand{\guG}{g_u^\Gamma}

\newcommand{\Rubar}{\,\underline{\!R\!}\,}

\newcommand{\tow}{\stackrel{\mathrm{w}}{\longrightarrow}}

\newcommand{\halmos}{$\sqcup\!\!\!\!\sqcap$}

\numberwithin{equation}{section}
\numberwithin{thm}{section}
\numberwithin{lemma}{section}
\numberwithin{prop}{section}
\numberwithin{remark}{section}
\numberwithin{cor}{section}
\numberwithin{ex}{section}

\begin{document}

\date{}
\title{General tax structures for a L\'evy  insurance risk process under the Cram\'er condition.}
\author{Philip S. Griffin\footnote{This work was partially supported by Simons Foundation Grant 226863.}\\
Syracuse University
}
\maketitle

\begin{abstract}
We investigate the L\'evy insurance risk model with tax under Cram\'er's condition.  A direct analogue of Cram\'er's estimate  for the probability of ruin in this model is obtained, together with the asymptotic distribution, conditional on ruin occurring, of several variables of interest related to ruin including the surplus immediately prior to ruin (undershoot) and shortfall at ruin (overshoot).  We also compute the present value of all tax paid conditional on ruin occurring.  The proof involves first transferring results from the model with no tax to the reflected process, and from there to the model with tax.  In doing so we also derive new results for the reflected process. 
\end{abstract}




\vspace*{10pt} \setcounter{equation}{0}
\section{Introduction}\label{s1}

Let $X$ be a L\'evy process which satisfies Cram\'er's condition,
\be\label{C}
Ee^{\ga X_1}=1\ \text{and} \  EX_1e^{\ga X_1}<\infty\text{ for some } \ga>0.
\ee
Such processes are widely used to model insurance risk.  In this setting $X$ represents the claims surplus process which is the excess in claims over premium.  The insurance company starts with initial capital $u$ and ruin occurs if $X$ exceeds $u$.  
The literature contains many papers considering various aspects of the behavior of $X$ at and approaching ruin.  
This paper investigates analogous questions in the presence of tax.

Let $\Gamma$ be a measurable process taking values in $[0,1]$ which is adapted to the usual augmented filtration of $X$ and set
\be\label{RG}
\RG_t=X_t-\int_0^t\Gamma_s d\Xubar_s=X_t+\int_0^t\Gamma_s d|\Xubar|_s,\ \ t\ge 0,
\ee
where $\Xubar_t=\inf_{s\le t}X_s$.  
If  $\Gamma\equiv 1$, then the resulting process is the reflected process of $X,$ reflected in its infimum, and in this case we write $R$ for $\RG$.
In general $\RG$ models a loss-carried-forward tax scheme in which the insurance company only pays tax when it is making a profit, by which we mean $\Xubar$ is decreasing, the amount of tax paid up to time $t$ being $\int_0^t\Gamma_s d|\Xubar|_s.$\footnote{Alternatively, one can view the payment as a profit participation scheme in which the cumulative payout up to time $t$ is $\int_0^t\Gamma_s d|\Xubar|_s.$} In the special case $\Gamma\equiv 0$ no tax is paid and $\RG_t=X_t,$  while if $\Gamma\equiv 1$ then all profits are paid and $\RG_t=R_t.$  When $\Gamma$ is a constant in $(0,1)$, a fixed proportion of the profits are paid as tax, and this case was considered by Albrecher and Hipp \cite{AH} and by Albrecher, Renaud and Zhou \cite{ARZ}.  A more general tax structure of the form $\Gamma_s=\gamma(\Xubar_s),$ was subsequently studied by Kyprianou and Zhou \cite{KZ}; see also Albrecher, Borst, Boxma and Resing \cite{ABBR} and Renaud \cite{R}.  In each of these papers $X$ is assumed to be spectrally positive and emphasis is placed on finding exact formulas for fixed $u$.  In the spirit of Cram\'er's original estimate, the aim here is to investigate asymptotic formulas under \eqref{C} without need of spectral positivity.

To describe the results we need to introduce a little notation.  For 
$u>0$ set 
\be\label{Ggt}
\tuG=\inf\{t:\RG_t>u\}.
\ee
If $\Gamma\equiv 0$ we simply write $\tuG=\tu$ and if $\Gamma\equiv 1$ we write $\tuG=\tu^R$.  We will usually refer to $\tuG$ as the ruin time.  
Here and throughout the paper we assume: 
\vskip.05in
\noindent {\bf Assumption C}:  \eqref{C} holds and $X$ is non-lattice.
\vskip.05in
\noindent In the lattice case  limits need to be taken through lattice values. The first result is an analogue of Cram\'er's estimate in this model. 

\begin{thm}\label{1tugn}  For some known constant $ \Upsilon_\Gamma>0,$
\be\label{1tuglim}
\lim_{u\to\infty} e^{\ga u}{P(\tuG<\infty)}= \Upsilon_\Gamma.
\ee
A sufficient condition for $\Upsilon_\Gamma$ to be finite is that $\Gamma$ be bounded away from $1$.
\end{thm}

When $\Gamma\equiv 0$ this reduces to Cram\'er's famous estimate.  The precise value of 
$\Upsilon_\Gamma$  is given in Theorem \ref{tugn}.  It is clear that some restriction is required on $\Gamma$ to ensure that $\Upsilon_\Gamma<\infty$ since when $\Gamma\equiv 1$  we are dealing with the reflected process in which case $P(\tu^R<\infty)=1$ for all $u$.

We next turn to the behavior of $\RG$ at and approaching ruin.  Outside of the case where $\RG=X$ the only known result is when $\RG=R$.  In this case Mijatovic and Pistorius \cite{MP} derived the joint limiting distribution of the undershoot and overshoot of the reflected process over level $u$ as $u\to\infty$. One interesting aspect of their result is that the joint limiting distribution obtained is the same as the corresponding limit for $X$ conditioned to cross over arbitrarily high levels.  To be precise they showed
\be\label{MP}
\lim_{u\to\infty}P(u-R_{\tu^R-}\in dy, R_{\tu^R}-u\in dx)=\lim_{u\to\infty}\Pu(u-X_{\tu-}\in dy,X_{\tu}-u\in dx)
\ee
in the sense of weak convergence, where $\Pu(\cdot)=P(\cdot|\tu<\infty).$  They also gave an expression for the limiting distribution.  Roughly speaking their approach was to prove that the limiting distribution on the right side exists and evaluate it, and then use this to show the limit on the left side exists and is the same.  We take this idea in the opposite direction and use the existence of the limit for the reflected process to
show that  if $\Gamma$ is bounded away from $1,$ then
\be\label{MPG}
\lim_{u\to\infty}\PuG(u-\RG_{\tuG-}\in dy, \RG_{\tuG}-u\in dx)=\lim_{u\to\infty}P(u-R_{\tu^R-}\in dy, R_{\tu^R}-u\in dx)
\ee
where $\PuG(\cdot)=P(\cdot|\tuG<\infty).$

In fact we are able to do much more.  Given a path $w$ in the usual Skorohod space $D$, let
\be
\tu=\inf\{t:w_t>u\} 
\ee 
and
\be
g_u=\sup\{s<\tu:w_s=\wbar_s\}
\ee
where $\wbar_t=\sup_{s\le t}w_s.$  Define a new path $w^u$ by
\be\label{Xu}
w^u_t= w_{(g_u+t)\wedge \tu} - \wbar_{g_u}, \quad t\ge 0.
\ee 
Thus, if $g_u<\tu<\infty,$ 
$w^u$ is the excursion of $w$ below its supremum  which results in $w$ first crossing over level $u$, and  $u-\wbar_{g_u}$ is the distance below level $u$ from which this excursion emanates. To simplify notation we will write  $\guG$ for $g_u(\RG)$, $\tuG$ for $\tu(\RG)$ (which agrees with \eqref{Ggt}) and 
$\RGu$ for $(\RG)^u$.  Let $\geps$ be the excursion of $R$ away from $0$ which first crosses over level $u$ and let $g^\geps_u=g_u(\geps)$ and
$
\tu^\geps=\tu(\geps).
$
Of course $\geps$ itself depends on $u$ but to simplify notation we suppress this dependence.
We then have

\begin{thm}\label{1ThTS} For suitable functions $G:[0,\infty)\times D\to [0,\infty)$, if
$\Gamma$ is bounded away from 1,  then
\be\ba\label{1nuPu}
\lim_{u\to\infty}\EuG G(u-\RbarG_{\guG}, \RGu)=\lim_{u\to\infty}EG(u-\gepsbar_{g^\geps_u}, \geps^u)
\ea\ee
\end{thm}

The starting point for the proof of Theorem \ref{1ThTS} is a result taken from \cite{G1}, see Theorem \ref{ThmGC} below.  It states the existence of the limit on the left side of \eqref{1nuPu} when $\RG=X,$ gives precise conditions on the function $G$ under which the limit exists and finally evaluates the limit in terms of the excursion measure of $R.$  In Theorem \ref{ThmC} this result is used to prove that the right side of \eqref{1nuPu} exists and the limit is the same as the limit from Theorem \ref{ThmGC}.  The result for $\geps$ is then used to prove the result in Theorem \ref{1ThTS} for general $\RG$.

As an example take 
$
G(y,w)=f(y-w_{\gz-},w_{\gz}-y)I(\gz<\infty)
$
where $f$ is bounded and continuous and $\gz=\inf\{t:w_s=w_t \text{ for all } s\ge t\}$ is the lifetime of $w.$
Then
\be
G(u-\wbar_{g_u}, w^u)=f(u-w_{\tu-},w_{\tu}-u)
\ee 
if $\tu<\infty,$ and so  \eqref{1nuPu} just reduces to \eqref{MPG} since the undershoot and overshoot of $u$ are the same for  $R$ as they are for $\geps.$  More generally, Theorem \ref{1ThTS} shows that the behavior of $\RG$ as ruin approaches, that is, over the time interval $[\guG,\tuG],$ is the same asymptotically as the behavior of $\geps$ over the interval $[g_u^\geps, \tu^\geps].$  Theorem \ref{1ThTS} would not be true if $\geps$ were replaced by $R$ since $g_u^R$ need not occur on the excursion $\geps.$  Combined  with the results in \cite{G1}, this allows us to derive the joint asymptotic distribution of several variables related to ruin; see Proposition \ref{P41}.  We can also evaluate some unbounded functionals of the path, for example certain Gerber-Shiu expected discounted penalty functions; see Proposition \ref{P42}.

Our final result concerns the limit of the expected present value of all tax payments conditional on ruin occurring.  In this result $\gd$ represents the discount rate and $\whk$ is the Laplace exponent of the descending bivariate ladder height process; see \eqref{kapdef}.
\begin{thm} 
If $\Gamma\equiv \gamma\in (0,1)$ and $X$ is spectrally positive, then for any $\gd\ge0$,
\be\label{div13}
\lim_{u\to\infty} \EuG \left(\int_0^{\tuG} e^{-\gd s}\Gamma_s d|\Xubar|_s\right)
= \frac{\gamma}{\ga(1-\gamma)+\whk(\gd,0)}.
\ee
\end{thm}  
When $\gd=0$ this gives the limit of expected total amount of tax paid conditional on ruin occurring.  In Theorem \ref{Div}, 
we derive the existence of a finite limit in \eqref{div13} when spectrally positivity is dropped and $\Gamma$ is only assumed to be bounded away from $1.$


 \setcounter{equation}{0}
\section{ Preliminaries }\label{s2}

Let $X=\{X_{t}: t \geq 0 \}$,
$X_0=0$,  be a L\'{e}vy process
with characteristic function $Ee^{i\theta X_{t}} = e^{t \Psi_X(\theta)}$,
given by the L\'{e}vy-Khintchine
representation
\ben\label{lrep}
\Psi_X(\theta) =
 \rmi\theta \gamma - \sigma^2\theta^2/2+
\int_{\R}(e^{\rmi\theta x}-1-
\rmi\theta x \mathbf{1}_{\{|x|<1\}})\Pi_X(d x),
\ {\rm for}\  \theta \in \mathbb{R}
\een
where $\Pi_X$ is the L\'evy measure of $X.$

Let  $D$  be the Skorohod space of functions  $w:[0,\infty )\to\R$ which are
right continuous with left limits.  
The usual right continuous completion of the filtration generated by the coordinate maps will be denoted by $\{{\cal F}_t\}_{t\ge 0}$.
$P_z$ denotes the probability measure induced on
${\cal F}=\vee_{t\ge 0} {\cal F}_t$
by the L\'evy process starting at $z\in\R$.
We usually write $P$ for $P_0$.
Set
\be\label{wwb}
\wbar_t=\sup_{s\le t} w_s,\quad \wbar=\sup_{s\ge 0} w_s
\ee
and 
\be\label{tu}
\tau_z=\inf\{t>0: w_t>z\}, \qquad   \tau^-_z=\inf\{t>0: w_t<z\}. 
\ee
It is convenient to assume that $X$ is given as the coordinate process on $D,$ and we will write
$X$ or $w$ depending on which seems clearer in the context.  This means for example that 
$\tz=\inf\{t:X_t>z\},$ in agreement with the notation in the introduction.

Let $(L^{-1}_s,H_s)_{s \geq 0}$
denote the 
ascending bivariate ladder
process of $X$, see Chapter VI of \cite{bert}.
When $X_t\to -\infty$ a.s., as is the case when \eqref{C} holds, $(L^{-1},H)$ is defective and may
be obtained from  a non-defective process
$({\cal L}^{-1},{\cal H})$ by exponential killing at some 
rate $q > 0$.
The bivariate renewal function of
$(L^{-1},H)$ is
 \begin{equation}\label{Vkdef}
V(t,x)= \int_0^\infty e^{-qs}P({\cal L}^{-1}_s\le t,{\cal H}_s\le x)ds,
\end{equation}
and its Laplace exponent $\kappa$ is given by
\begin{equation} \label{kapdef}
e^{-\kappa(a, b)}  =  e^{-q}E e^{ -a{\cal L}^{-1}_1-b{\cal H}_1}
 \end{equation}
for values of $a,b\in \R$ for which the expectation is finite.

Let $\whX_t=-X_t$, $t\ge 0$ denote the dual process, and $(\whL^{-1}, \whH)$ the corresponding  ascending bivariate
ladder process of $\whX$.  Its  bivariate renewal function will be denoted by $\whV$, and its Laplace exponent by
$\widehat \kappa.$ 
We assume the local times $L$ and $\whL$ are chosen so that
the Wiener-Hopf factorization takes the form
 \begin{equation}\label{WH}
{\kappa(a,-ib) \wh\kappa(a,ib)
= a-\Psi_X( b),\  a\ge 0,b\in\R.} 
  \end{equation}

Recall
\ben
\Xubar_t=\inf_{s\le t}X_s
\een
and $R=X-\Xubar$ is the reflected process of $X$.  Its local time at $0$ is $\whL$.  Define $e_t$ by 
\ben
e_t(s)=R_{(\whL_{t-}^{-1}+s)\wedge\whL_{t}^{-1}}=\left(X_{(\whL_{t-}^{-1}+s)\wedge\whL_{t}^{-1}}-\Xubar_{\whL_{t-}^{-1}}\right)\vee 0.
\een
Then $e_t$ is the excursion of $R$ away from $0$ at local time $t,$ and
$\{(t,e_t):\whL^{-1}_{t-}<\whL^{-1}_{t}\}$ is an ${\cal F}_{\whL_t^{-1}}-$Poisson point process with intensity measure $\whn$, called the excursion measure.  This measure is defined on ${\cal F}$ but is supported on the set
\be
\{w\in D:w_s\ge 0 \text{ for all } s\}.
\ee
For $t\ge 0,$ in agreement with the notation in \eqref{wwb}, we write
\ben
\ebar_t(r)
=\sup_{s\le r} e_t(s)\ \ \text{and }\   \ebar_t=\sup_{s\ge 0} e_t(s).
\een
Thus $\ebar_t$ is the height of the excursion at local time $t$. 
For $u>0$ let
\ben
T_u=\inf\{t:\ebar_t>u\}.
\een
Thus $e_{T_u}$ is the first excursion to cross over level $u$.  This is the excursion $\geps$ referred to in the introduction. 
For convenience we will use both notations $e_{T_u}$ and $\geps$ going forward.


 \setcounter{equation}{0}
\section{The Reflected Process}\label{sC}

We briefly review some well known properties of L\'evy processes which satisfy \eqref{C}.  Since $e^{\ga X_t}$ is a non-negative martingale, it converges, and so $X_t\to-\infty$ a.s.  Define a new probability measure $Q$ on ${\cal F}$ by 
\be\label{Q}
dQ=e^{\ga X_t} dP \ \text{ on } {\cal F}_t.
\ee
Under $Q$, $X_1$ has positive mean and so $X_t\to\infty$ a.s.
For any stopping time $\tau$, \eqref{Q}
continues to hold if $t$ is replaced by  $\tau$ provided we restrict to $\{\tau<\infty\}.$  In particular for any $u>0$
\be\label{E=1}
E(e^{\ga X_{\tu}};\tu<\infty)=Q(\tu<\infty)=1.
\ee
Thus
\be\label{ebd}
E^{(u)}e^{\ga(X_{\tu}-u)}=\frac{e^{-\ga u}E{(e^{\ga X_{\tu}};\tu<\infty)}}{P(\tu<\infty)}=\frac 1{e^{\ga u}P(\tu<\infty)}.
\ee
We will need the following results in which $\Upsilon\in(0,\infty)$ is a known constant; 
\begin{subequations}
\begin{align}
 \lim_{u\to\infty}e^{\ga u}P(\tu<\infty)&= \Upsilon, \label{MP1}\\
 \lim_{u\to\infty}e^{\ga u}\whn(\tu<\gz)&=\Upsilon \whk(0,\ga),\label{MP3}\\
 \lim_{u\to\infty}\whn(e^{\ga\gep_{\tu}}; \tu<\gz)&=\whk(0,\ga),\label{MP5}\\
\end{align}
\end{subequations}
where, as is customary,  we write  $\whn(K;B)$ for $\int_BK(w)\whn(dw).$
Equations \eqref{MP1}, \eqref{MP3}  and  \eqref{MP5} are (5) in \cite{BD94}, Theorem 1  in \cite{DM}  and (3.30) in \cite {MP} respectively.   Two simple consequences of these results are
\be\ba\label{Pxn}
\sup_{u>0}\frac{P(\tau_{u-x}<\infty)}{\whn(\tu<\gz)}&\le Ce^{\ga x} 
\ea\ee
for every $x\ge 0$ 
where $C$ does not depend on $x$, and 
\be\label{whnbd}
\whn(e^{\ga\gep_{\tv}}; \tv<\gz)<\infty \text{ for all } v>0.
\ee
We will use $C$ throughout to denote the value of an absolute constant whose precise value is not important and may change from one usage to the next.

The excursion measure $\whn$ of $X-\Xubar$ was introduced in Section \ref{s2}.  To obtain results about $\geps^u=e_{T_u}^u$, the final excursion of $\geps$ below its supremum prior to first passage over $u$, we will also need the excursion measure  
$n$ of $X-\Xbar.$ This is defined similarly to $\whn$, but to compute the limiting overshoot, we need to include additional information about any possible jump at the end of the excursion. For this reason we let
\ben
\whe_t(s)=X_{(L_{t-}^{-1}+s)\wedge L_{t}^{-1}}-\Xbar_{L_{t-}^{-1}}.
\een
Then  $\{(t,\whe_t):L^{-1}_{t-}<L^{-1}_{t}\}$ is a Poisson point process with intensity measure $n.$
To account for the possibility that first passage of $X$ over a level is due to a jump which does not occur at the end of an excursion interval we need to add an additional term to $n$.  For this  let ${\bf x}\in D$ be the path which is identically $x\in \R$ and define
\be
\nbar(A)=n(A)+d_{L^{-1}}\Pi_X^+(\{x:{\bf x}\in A\}),\  \text{for } A\in {\cal F},
\ee 
where $d_{L^{-1}}$ is the drift of $L^{-1}$ and $\Pi_X^+(\cdot)=\Pi_X(\cdot\cap(0,\infty))$.  The following result,  Theorem 3.1 in \cite{G1}\footnote{
In \cite{G1}, $(\whL^{-1}, \whH)$ is taken to be the weakly ascending bivariate ladder height process of $\whX$.  
This distinction is only an issue when $X$ is compound Poisson, in which case  a slightly different definition of $X^u$ is used in \cite{G1}, but the proofs given here will still work with minor modifications in that case.}, indicates the need to introduce $\nbar$ and is the starting point for our investigation.

\begin{thm}\label{ThmGC} Assume $G:[0,\infty)\times D\to [0,\infty)$ is measurable, $G(y, w)e^{-\ga(w_\gz-y)}I(w_\gz>y)$ is bounded in $(y,w)$, and $\nbar(B_y^c)=0$ for a.e. y  with respect to Lebesgue measure, where
\ben
B_y=\{w: G(\cdot\ ,w) \text{ is continuous at } y\}.
\een  
Then, with $d_{H}$ denoting the drift of $H$, we have
\be\ba\label{lim31}
\lim_{u\to\infty}E^{(u)}G(u-\Xbar_{g_u}, X^u)
&=  \int_{[0,\infty) }\frac{\ga}q e^{\ga y}dy \int_{{D}}G(y,\gep)\ \nbar(d\gep, \gep_\gz>y)+d_H\frac{\ga}q G(0,{\bf 0}),
\ea\ee
where the limit is finite.
\end{thm} 

The final term in the limit accounts for the possibility of $X$ creeping over level $u$.
The continuity condition holds when $G$ is continuous in $y$ for each $w$.  The boundedness condition holds when $G$ is bounded, but allows for certain unbounded functions.  Since $w^u_\gz=w_{\tu}-\wbar_{g_u}$, the boundedness condition also implies
\be\label{Gbd}
G(u-\wbar_{g_u}, w^u)\le Ce^{\ga(w_{\tu}-u)}\ \text {on } \{\tu<\infty\},
\ee
and hence it follows immediately from \eqref{ebd} and \eqref{MP1} that
\be\label{Gbd1}
\sup_{u>0} E^{(u)}G(u-\Xbar_{g_u}, X^u)<\infty.
\ee

We are now ready to state the main result of this section.
In its proof we will make use of the following well known property of Poisson point processes, see \cite{bert};  for any nonnegative measurable functional $K:D\to\R$
\be\label{pppu}
EK(\geps)=EK( e_{T_u})=\frac{\whn(K;  \tu<\gz)}{\whn(\tu<\gz)}.
\ee
Additionally we make use of ideas from \cite{MP},  together with Theorem 
\ref{ThmGC}.

\begin{thm}\label{ThmC} Assume $G$ satisfies the same conditions as in Theorem \ref{ThmGC}, then 
\be\ba\label{nuPu32}
\lim_{u\to\infty}EG(u-\gepsbar_{g^\geps_u}, \geps^u)
&=\lim_{u\to\infty}E^{(u)}G(u-\Xbar_{g_u}, X^u). 
\ea\ee
\end{thm}

\noindent{\bf Proof}\ \ 
For any $0<v\le u<\infty$, by \eqref{pppu} and the 
strong Markov property for $\whn$,   
\be\ba
EG(u-\gepsbar_{g^\geps_u}, \geps^u)
&=\frac{\whn(G(u-\wbar_{g_u}, w^u);  \tu<\gz)}{\whn(\tu<\gz)}\\
&=\int_{[v,u]}\whn(\gep_{\tv}\in dx,\tv<\zeta)\frac{E_x(G(u-\Xbar_{g_u}, X^u); \tu<\too^-)}{\whn(\tu<\gz)}\\
&\qquad +\frac{\whn(G(u-\wbar_{g_u}, w^u); \gep_{\tv}>u, \tu<\gz)}{\whn(\tu<\gz)}\\
&=I+II.
\ea
\ee

For $II$, by \eqref{Gbd} we have
\be\ba
II&\le \frac{C\whn(e^{\ga(w_{\tu}-u)};\gep_{\tv}>u,\tu<\zeta)}{\whn(\tu<\zeta)}= \frac{C\whn(e^{\ga w_{\tv}};\gep_{\tv}>u,\tv<\zeta)}{e^{\ga u}\whn(\tu<\zeta)}
\ea\ee
since $\tu=\tv$ if $w_{\tv}>u$.  Thus for every $v>0$, by \eqref{MP3}, \eqref{whnbd} and dominated convergence
\be\label{3II}
\lim_{u\to\infty} II=0.
\ee

For $I$, we begin by writing
\be\ba
I&=\int_{[v,u]}\whn(\gep_{\tv}\in dx,\tv<\zeta)\frac{E_x(G(u-\Xbar_{g_u}, X^u); \tu<\infty)}{\whn(\tu<\gz)}\\
&\qquad -\int_{[v,u]}\whn(\gep_{\tv}\in dx,\tv<\zeta)\frac{E_x(G(u-\Xbar_{g_u}, X^u); \too^-<\tu<\infty)}{\whn(\tu<\gz)}\\
&=I_1-I_2.
\ea\ee
To find the limit of $I_1$ we first write
\be\ba
I_1&=\int_{[v,u]}\whn(\gep_{\tv}\in dx,\tv<\zeta)\frac{P( \tau_{u-x}<\infty)}{\whn(\tu<\gz)}f(u-x)
\ea\ee
where
\be
f(u)=\frac{E(G(u-\Xbar_{g_u}, X^u); \tu<\infty)}{P( \tu<\infty)}=E^{(u)}G(u-\Xbar_{g_u}, X^u).
\ee
By \eqref{Gbd1} and  Theorem \ref{ThmGC},  $f$ is bounded and $f(u)$ converges as $u\to\infty$ to the right hand side of \eqref{lim31}. 
By  \eqref{MP1} and \eqref{MP3}, for every $x\ge 0$
\be\ba
\frac{P(\tau_{u-x}<\infty)}{\whn(\tu<\gz)}&\to \frac{e^{\ga x}}{\whk(0,\ga)}
\ea\ee
as $u\to\infty.$  Thus by \eqref{Pxn} and \eqref{whnbd}, for every $v>0,$ we can apply dominated convergence to obtain
\be
\lim_{u\to\infty}I_1=
\frac{\whn(e^{\ga\gep_{\tv}}; \tv<\gz)}{\whk(0,\ga)}\left(\int_{[0,\infty) }\frac{\ga}q e^{\ga y}dy \int_{{D}}G(y,\gep)\ \nbar(d\gep, \gep(\gz)>y)+d_H\frac{\ga}q G(0,{\bf 0})\right).
\ee
Letting $v\to\infty$ and using \eqref{MP5} then gives
\be\label{3I1}
\lim_{v\to\infty}\lim_{u\to\infty}I_1=\int_{[0,\infty) }\frac{\ga}q e^{\ga y}dy \int_{{D}}G(y,\gep)\ \nbar(d\gep, \gep(\gz)>y)+d_H\frac{\ga}q G(0,{\bf 0}).
\ee
For $I_2$ we first observe that for any $x\ge 0,$ by the strong Markov property  and then \eqref{E=1}
\be\ba
{E_x (e^{\ga(X_{\tu}-u)}; \too^-<\tu<\infty)}&
=\int_{z\le 0} P_x(X_{\too^-}\in dz;\too^-<\tu)E_z(e^{\ga(X_{\tu}-u)}; \tu<\infty)\\
&=\int_{z\le 0} P_x(X_{\too^-}\in dz;\too^-<\tu) e^{-\ga(u-z)}\\
&\le e^{-\ga u}.
\ea\ee
Thus by \eqref{Gbd}
\be\ba
I_2&\le C\int_{[v,u]}\whn(\gep_{\tv}\in dx,\tv<\zeta)\frac{E_x (e^{\ga(X_{\tu}-u)}; \too^-<\tu<\infty)}{\whn(\tu<\gz)}\\
&\le \frac{C\whn(\tv<\zeta)}{e^{\ga u}\whn(\tu<\gz)}\\
&\to \frac{C\whn(\tv<\zeta)}{\Upsilon\whk(0,\ga)}
\ea\ee
as $u\to\infty$  by \eqref{MP3}.  Hence
\be
\lim_{v\to\infty}\lim_{u\to\infty}I_2=0.
\ee
Combined with \eqref{3II} and \eqref{3I1} this completes the proof.
\qquad\halmos
\vskip.1in

It is now immediate that any limit laws computed from $(u-\Xbar_{g_u}, X^u)$ under $\Pu$ translate to limit laws for 
$(u-\gepsbar_{g^\geps_u}, \geps^u)$ under $P$.  
As one example, from Theorem 3.4 of \cite{G1}, we have the following extension of the Theorem 2 in \cite{MP};

\begin{prop}\label{P31}  For any $y\ge0, x\ge 0, v\ge 0, t\ge 0$
\be\ba\label{jtconC}
P(&u-\gepsbar_{g^\geps_u}\in dy,\geps_{\tu^\geps}-u\in dx, u-\geps_{\tu^\geps-}\in dv, \tu^\geps- g^\geps_u\in dt)\\
&\tow\frac{\ga }q e^{\ga y} dy I(v\ge y) \whV(dt,dv-y)\Pi_X(v+dx) + d_H\frac{\ga }q \delta_{(0,0,0,0)}(dx,dy,dv,dt)
\ea\ee
where $\delta_{(0,0,0,0)}$ is a point mass at ${(0,0,0,0)}.$
\end{prop}

Integrating out $y$ and $t$  yields the joint limiting distribution of the undershoot and overshoot of $R$ obtained in \cite{MP}.  However, as pointed out in the introduction, we can not in general replace $\geps$ by $R$ in Proposition  \ref{P31} since $g_u^R$ need not occur on the excursion $\geps$.


 \setcounter{equation}{0}
\section{General Tax Structures}\label{sT}

From the definition of $\RG$ in \eqref{RG}, 
it is easily seen, see Lemma 2.1 of \cite{KZ}, that
\be\label{RGu}
\Rubar_t^{\Gamma}=\Xubar_t-\int_0^t\Gamma_sd\Xubar_s,
\ee 
where $\Rubar_t^{\Gamma}=\inf_{s\le t}\RG_s$.  Thus 
\be
R_t^{\Gamma}-\Rubar_t^{\Gamma}=R_t.
\ee
As a consequence the excursions of $R^\Gamma$ above its infimum do not depend of $\Gamma.$   They are precisely the same as the excursions of $R,$  a property which plays a central role in our use of excursion theory in this section.  
It is useful to picture $\RG$on the interval $(\whL^{-1}_{t-},\whL^{-1}_{t})$ as the excursion $e_t$  emanating from level
$\Rubar^{\Gamma}_{\whL^{-1}_{t-}},$  that is
\be\label{eem}
\RG_{\whL^{-1}_{t-}+s}=\Rubar^{\Gamma}_{\whL^{-1}_{t-}}+e_t(s),\ \ 0\le s<  \whL^{-1}_{t}-\whL^{-1}_{t-}.
\ee

Set
\be\label{HGdef}
\whH^{\Gamma}_t=-\Rubar^{\Gamma}_{\whL^{-1}_{t}}=-\int_0^{\whL^{-1}_{t}}(1-\Gamma_s)d\Xubar_s=\int_0^{\whL^{-1}_{t}}(1-\Gamma_s)d|\Xubar|_s.
\ee
This plays the same role in relation to $\RG$ that $\whH$ does in relation to $X.$   We will often use, without further mention, that for each fixed $t$, $\whH^{\Gamma}_t=\whH^{\Gamma}_{t-}$ a.s. This is because
\be\ba
0\le \whH^{\Gamma}_t-\whH^{\Gamma}_{t-}&= \lim_{r\uparrow t} \int_{({\whL^{-1}_{r}},{\whL^{-1}_{t}}]}(1-\Gamma_s)d|\Xubar|_s\\
&\le \lim_{r\uparrow t} (|\Xubar|_{\whL^{-1}_{t}}-|\Xubar|_{\whL^{-1}_{r}}) \\
&=\lim_{r\uparrow t} (\whH_t-\whH_r)\\
&=\whH_t-\whH_{t-}=0 \text{ a.s.}
\ea\ee

\begin{lemma}\label{mcT}
Let $\mcT=\{t:\whL^{-1}_{t-}< \whL^{-1}_{t}\}.$ 
Then $P(\whH^\Gamma_{t-}=-\Rubar^{\Gamma}_{\whL^{-1}_{t-}} \text{ for all $t\in \mcT$})=1.$
\end{lemma}

\noindent{\bf Proof}\ \ By a standard argument similar to p161 of \cite{bert},
\be\label{Xlcts}
P(\Xubar_{\whL^{-1}_{t-}} \text{ is left continuous at all } t\in\mcT)=1.
\ee 
Hence, on this set, by \eqref{RGu} and \eqref{HGdef}, for any $r<t$
\be\ba
0\le |\Rubar^{\Gamma}_{\whL^{-1}_{t-}}|-\whH^\Gamma_{t-}&\le |\Rubar^{\Gamma}_{\whL^{-1}_{t-}}|-\whH^\Gamma_{r}\\
&=\int_{({\whL^{-1}_{r}}, {\whL^{-1}_{t-}}]}(1-\Gamma_s)d|\Xubar|_s\\
&\le |\Xubar|_{\whL^{-1}_{t-}}-|\Xubar|_{\whL^{-1}_{r-}}\to 0 
\ea\ee
as $r\to t$ if $t\in \mcT.$
\qquad\halmos
\vskip.1in

Thus we may replace $\Rubar^{\Gamma}_{\whL^{-1}_{t-}}$ in \eqref{eem} with 
$-\whH^\Gamma_{t-}$. With this we are ready to obtain the analogue of Cram\'er's estimate for $\tuG.$

\begin{thm}\label{tugn}  
We have
\be\label{tuglim}
\lim_{u\to\infty} e^{\ga u}{P(\tuG<\infty)}= \Upsilon\whk(0,\ga)\int_0^\infty Ee^{-\ga \whH_t^\Gamma}dt.
\ee
A sufficient condition for the limit to be finite is that $\Gamma$ be bounded away from $1$.  If $\Gamma\equiv \gamma$ for some fixed 
$\gamma\in [0,1),$  then 
\be\label{Cramg}
e^{\ga u}P(\tuG<\infty)\to \frac{\Upsilon\whk(0, \ga)}{\whk(0, \ga(1-\gamma))}.
\ee
\end{thm}

\noindent{\bf Proof}\ \  
For notational convenience we set $\mcRG_{t}=\RbarG_{\whL^{-1}_{t-}} $  Note that 
$
\mcRG_{t-}=\RbarG_{\whL^{-1}_{t-}-}. 
$ 
Further $\mcRG_{t-}$ and $\whH^\Gamma_{t-}$ are both left continuous and ${\cal F}_{\whL_t^{-1}}-$adapted, hence for any $u>0$ by the compensation formula applied to the point process of excursions
\be\ba
P(\tuG<\infty)&= E \sum_tI(\mcRG_{t-}\le u,  
\ \ebar_t>u+\whH^\Gamma_{t-})\\
&=\int_0^\infty dt\int_z P(\mcRG_{t-}\le u, 
\whH^\Gamma_{t-}\in dz)\whn(\tau_{u+z}<\gz).
\ea\ee
By \eqref{MP3}, there is a constant $C$ such that for all $u\ge 1$ and all $z\ge 0$
\be\label{npot}
\frac{\whn(\tau_{u+z}<\gz)}{\whn(\tau_{u}<\gz)}\le Ce^{-\ga z}.
\ee
Thus if the integral in \eqref{tuglim} is finite then by dominated convergence and \eqref{MP3}  
\be\label{tugnlim}
\lim_{u\to\infty} \frac{P(\tuG<\infty)}{\whn(\tu<\gz)}= \int_0^\infty Ee^{-\ga \whH_t^\Gamma}dt.
\ee
If the integral is infinite then \eqref{tugnlim} follows from Fatou.  In either case \eqref{tuglim} then follows from \eqref{tugnlim} and \eqref{MP3}.

If $\Gamma$ is bounded away from 1, then  $\Gamma\le\gamma$ for some $ \gamma \in [0,1).$  Since, by \eqref{HGdef}, $\whH^\Gamma$ is decreasing in $\Gamma$, to prove the limit in \eqref{tuglim} is finite, it suffices to prove it is finite when $\Gamma\equiv \gamma.$
 Thus assume $\Gamma\equiv \gamma$ for some fixed 
$\gamma\in [0,1).$  Then by \eqref{HGdef},
$
\whH^{\Gamma}_t=(1-\gamma)\whH_t
$
and
\be\label{HGint}
\int_0^\infty Ee^{-\ga \whH_t^\Gamma}dt=\int_0^\infty e^{-\whk(0, \ga(1-\gamma))t}dt=\frac 1{\whk(0, \ga(1-\gamma))}
\ee
which also proves \eqref{Cramg}.

\qquad\halmos\vskip.1in

If $\Gamma\equiv 1$, then $\whH^\Gamma\equiv 0$ so the integral in \eqref{tuglim} is infinite.  Of course in this case $\RG=R$ and $P(\tuG<\infty)=1,$ and so the limit in \eqref{tuglim} must be infinite. This example illustrates that some condition is required on $\Gamma$ if the limit in \eqref{tuglim} is to be finite.  A  simple sufficient condition, as shown in Theorem \ref{tugn},   is that $\Gamma$  be bounded away from $1$.  The following example shows that this is not necessary.

\begin{ex}\label{eg1} Fix $\gb>0$ and let $\Gamma_s=f(\Xubar_s)$ where
\be
f(s)=
\begin{cases}
0, & s\le \gb\\
1-{\gb}{s^{-1}}, & s\ge \gb.
\end{cases}
\ee
 Assume $X$ is spectrally positive. Then we may take $\whL_t=|\Xubar|_t$ in which case $\whH_t=t$.  Thus by \eqref{HGdef}, after a change of variable,
\be\label{HG1}
\whH_t^{\Gamma}=\int_0^t (1-f(s))ds=
\begin{cases}
t, & t\le\gb\\
\gb+\gb\ln(t/\gb), & t>\gb.
\end{cases}
\ee
Hence
\be\label{HG2}
\int_0^\infty Ee^{-\ga \whH_t^\Gamma}dt=
\begin{cases}
\infty, & \gb\le \ga^{-1}\\
c_\Gamma, & \gb> \ga^{-1}.
\end{cases}
\ee
where 
\be
c_\Gamma=\frac{\ga\gb-1+e^{-\ga\gb}}{\ga(\ga\gb-1)}.
\ee
Since $\whk(0,\ga)=\ga,$ \eqref{tuglim} then becomes
\be
\lim_{u\to\infty} e^{\ga u}{P(\tuG<\infty)}= \frac{\Upsilon(\ga\gb-1+e^{-\ga\gb})}{\ga\gb-1}
\ee
if $\ga\gb>1$ and infinite otherwise.

\end{ex}

We now turn to the proof of Theorem \ref{1ThTS} for which we first give a precise statement.

\begin{thm}\label{ThTS} Assume $G$ satisfies the same conditions as in Theorem \ref{ThmGC} and $\Gamma$ is bounded away from 1.  Then
\be\ba\label{nuPu}
\lim_{u\to\infty}\EuG G(u-\RbarG_{\guG}, \RGu)=\lim_{u\to\infty}EG(u-\gepsbar_{g^\geps_u}, \geps^u)
\ea\ee
\end{thm}

The proof of the Theorem will be accomplished via a series of propositions 
 which break down the expectations depending on whether $\tuG$ occurs during the first excursion  $e_{T_u}$
to reach level $u$ or on some later excursion, which necessarily must  also reach level $u$.
If $\tuG$ occurs on  the first excursion 
to exceed $u$, then $\whL^{-1}_{T_u-}\le \tuG< \whL^{-1}_{T_u}.$  
Set
\be
\suG=\inf\{t\ge \whL^{-1}_{T_u}:\RG_t>u\}.
\ee
The following result asserts that if $R^\Gamma$ ever exceeds $u$, it does so only on the first excursion 
that exceeds $u$, except on a set which is negligible relative to the size of the event $\{\tuG<\infty\}.$

\begin{prop}\label{FE}  Assume $\Gamma$ is bounded away from 1, then
\be\label{neg}
\lim_{u\to\infty} \frac{P(\suG<\infty)}{P(\tuG<\infty)}=0.
\ee
Additionally, for any $x\ge 0,$
\be\label{all}
\lim_{u\to\infty} \frac{P(\whL^{-1}_{T_u-}\le \tau_{u-x}^\Gamma\le \tuG< \whL^{-1}_{T_u})}{P(\tuG<\infty)}=1.
\ee
\end{prop}

\noindent{\bf Proof}\ \  
Choose $\gamma\in [0,1)$ such that $\Gamma\le \gamma$.  For notational convenience we will write $\whH^{\gamma}$ for the process $\whH^{\Gamma}$ when $\Gamma\equiv \gamma.$ Let
$\whK_u$ be independent of $X$ and have the same distribution as $\whH^{\gamma}_{T_u}=(1-\gamma)\whH_{T_u}.$  
The Poisson point process $\{(t,e_t):\whL^{-1}_{(T_u+t)-}<\whL^{-1}_{(T_u+t)}, t>0\}$ is independent of 
$\whH_{T_u}$
and 
its joint distribution with the increment process $\whH_{T_u+\cdot}-\whH_{T_u}$ is the same as the joint distribution of $\{(t,e_t):\whL^{-1}_{t-}<\whL^{-1}_{t}\}$ with $\whH.$
Thus, since $\whH^\Gamma$ is decreasing in $\Gamma,$
\be\ba
P(\suG<\infty)&= P(\ebar_{T_u+t}>u+\whH^\Gamma_{(T_u+t)-} \text{ for some } t>0)\\
&\le P(\ebar_{T_u+t}>u+\whH^\gamma_{(T_u+t)-} \text{ for some } t>0)\\
&= P(\ebar_t>u+(1-\gamma)\whH_{t-}+\whK_u\text{ for some } t> 0)\\
&=\int_zP(\ebar_t>u+(1-\gamma)\whH_{t-}+z \text{ for some } t> 0)P(\whK_u\in dz)\\
&\le \int_zE\sum_tI(\ebar_t>u+(1-\gamma)\whH_{t-}+z)P(\whK_u\in dz)\\
&= \int_z E\int_0^\infty dt\ \whn(\tau_{u+(1-\gamma)\whH_{t-}+z}<\gz)P(\whK_u\in dz).
\ea\ee
Thus by \eqref{npot}
\be\ba
\limsup_{u\to\infty} \frac{P(\suG<\infty)}{\whn(\tu<\gz)} 
&\le \left({C}\int_0^\infty Ee^{-\ga(1-\gamma) \whH_t}dt\right) \limsup_{u\to\infty} \int_0^\infty e^{-\ga z}P(\whK_u\in dz)\\
&= \frac{C}{\whk(0,(1-\gamma)\ga)} \limsup_{u\to\infty} \int_0^\infty e^{-\ga z}P((1-\gamma)\whH_{T_u}\in dz)\\
&=0
\ea\ee
since $T_u\to\infty$ a.s.  Together with \eqref{MP3} and Theorem \ref{FE} this proves \eqref{neg}.

For \eqref{all}, observe that 
\be\ba
\{\tuG<\infty\}&=\{\whL^{-1}_{T_u-}\le \tuG<\whL^{-1}_{T_u}\}\cup
\{\whL^{-1}_{T_u}\le \tuG<\infty\}\\
&\subseteq\{\whL^{-1}_{T_u-}\le \tau_{u-x}^\Gamma\le \tuG< \whL^{-1}_{T_u}\}
\cup\{\sigma_{u-x}^\Gamma<\infty \}\cup \{\whL^{-1}_{T_u}\le \tuG<\infty\},
\ea\ee
since if $\tau_{u-x}^\Gamma<\whL^{-1}_{T_u-}$ and $\tuG<\whL^{-1}_{T_u}$ then $\sigma_{u-x}^\Gamma<\infty.$
But
\be
\frac{P(\whL^{-1}_{T_u}\le \tuG<\infty)}{P(\tuG<\infty)}\le \frac{P(\sigma_{u}^\Gamma<\infty)}{P(\tuG<\infty)}\to 0
\ee
by \eqref{neg} and
\be
\frac{P(\sigma_{u-x}^\Gamma<\infty)}{P(\tuG<\infty)}
=\frac{P(\sigma_{u-x}^\Gamma<\infty)}{P(\tau_{u-x}^\Gamma<\infty)}
\frac{P(\tau_{u-x}^\Gamma<\infty)}{P(\tuG<\infty)}\to 0
\ee
by Theorem \ref{FE} and \eqref{neg}.  Hence \eqref{all} holds.
\qquad\halmos
\vskip.2in

 Let  $T_u^1=T_u$ and for $k\ge 2,$ $T_u^k=\inf\{t>T_u^{k-1}:\ebar_t>u\}.$  Then 
\be
\{ \tuG<\infty\}=\bigcup_{k=1}^\infty \{\whL^{-1}_{T^k_{u}-}\le \tuG< \whL^{-1}_{T^k_u}\}.
\ee
Further
\be\label{Tu1}
\{\whL^{-1}_{T_u-}\le \tuG< \whL^{-1}_{T_u}\}=\{\ebar_{T_u}>u+\whH^\Gamma_{T_u-}\},
\ee
and for $k\ge 2$,
\be\label{Tuk}
 \{\whL^{-1}_{T^k_{u}-}\le \tuG< \whL^{-1}_{T^k_u}\}=\{\mcRG_{T^k_u-}\le u, \ebar_{T^k_u}>u+\whH^\Gamma_{T^k_u-}\},
 \ee
where recall $\mcRG_{t}=\RbarG_{\whL^{-1}_{t-}}.$  
Set
\be
O_u=O_u(w)=w_{\tu} -u,
\ee
the overshoot of $u$ by the path $w$ on $\{\tu<\infty\}.$  Then 
by \eqref{Gbd}
\be\label{GO}
G(u-\RbarG_{\guG}, \RGu)\le C e^{\ga O_u(\RG)}\ \text{ if $\tuG<\infty,$ }
\ee
while for any $k\ge 1$
\be\label{os}
O_u(\RG)=O_{u+\whH^\Gamma_{T^k_u-}}(e_{T^k_u})\ \text{ if $ \whL^{-1}_{T^k_{u}-}\le \tuG< \whL^{-1}_{T^k_u}$.}
\ee

If $G$ were bounded the following result would follow immediately from Proposition \ref{FE}.  In place of boundedness we make use of the weaker condition \eqref{Gbd} that $G$ satisfies.

\begin{prop}\label{prop2} Assume $G$ satisfies the same conditions as in Theorem \ref{ThmGC} and $\Gamma$ is bounded away from 1.  Then
\be\ba\label{nuPu42}
\lim_{u\to\infty}\EuG[G(u-\RbarG_{\guG}, \RGu):\whL^{-1}_{T_u}\le \tuG<\infty]=0.
\ea\ee
\end{prop}

\noindent{\bf Proof}\ \  
Fix $k\ge 2$.  For any non-negative measurable function $f,$ by \eqref{Tuk}, \eqref{os} and the compensation formula
\be\ba\label{fO}
E[f(O_u&(\RG));\whL^{-1}_{T^k_{u}-}\le \tuG< \whL^{-1}_{T^k_u}]\\
&=E[f(O_{u+\whH^\Gamma_{T^k_u-}}(e_{T^k_u}));\mcRG_{T^k_u-}\le u, \ebar_{T^k_u}>u+\whH^\Gamma_{T^k_u-}]\\
&=E\sum_t I\left(\sum_{s<t}I(\ebar_s>u)=k-1, \mcRG_{t-}\le u, \ebar_t>u+\whH^\Gamma_{t-}\right)f(O_{u+\whH^\Gamma_{t-}}(e_t))\\
&=\int_0^\infty dt\int_z P\left(\sum_{s<t}I(\ebar_s>u)=k-1,\mcRG_{t-}\le u, \whH^\Gamma_{t-}\in dz\right)\whn(f(O_{u+z}); \tau_{u+z}<\gz).
\ea\ee
Setting $f\equiv 1$ gives
\be\label{PTk}
P(\whL^{-1}_{T^k_{u}-}\le \tuG<\whL^{-1}_{T^k_u})= \int_0^\infty dt\int_z P\left(\sum_{s<t}I(\ebar_s>u)=k-1,\mcRG_{t-}\le u, \whH^\Gamma_{t-}\in dz\right)\whn(\tau_{u+z}<\gz).
\ee
Now by \eqref{MP3} and \eqref{MP5}, for some $u_0$
\be\label{35}
\sup_{u\ge u_0}\frac{\whn(e^{\ga O_{u}}; \tau_{u}<\gz)}{\whn(\tau_{u}<\gz)}\le 2\Upsilon^{-1}
\ee
Thus taking $f(x)=e^{\ga x}$  and $u\ge u_0$ in \eqref{fO}, and using \eqref{PTk} and \eqref{35}, we have for any $k\ge 2$
\be
E[e^{\ga O_u(\RG)};\whL^{-1}_{T^k_{u}-}\le \tuG< \whL^{-1}_{T^k_u}]
\le 2\Upsilon^{-1} P(\whL^{-1}_{T^k_{u}-}\le \tuG< \whL^{-1}_{T^k_u}).
\ee
Hence by \eqref{GO}, for $u\ge u_0$
\be\ba
\EuG[G(u-\RbarG_{\guG}, \RGu);\whL^{-1}_{T_u}\le \tuG<\infty]
&\le C  \sum_{k=2}^\infty \EuG[e^{\ga O_u(\RG)};\whL^{-1}_{T^k_{u}-}\le \tuG<\whL^{-1}_{T^k_u}]\\
&\le 2\Upsilon^{-1}C  \sum_{k=2}^\infty \frac{P(\whL^{-1}_{T^k_{u}-}\le \tuG< \whL^{-1}_{T^k_u})}{P(\tuG<\infty)}\\
&= 2\Upsilon^{-1}C
\frac{P(\whL^{-1}_{T_u}\le \tuG<\infty)}{P(\tuG<\infty)}\to 0
\ea\ee
by \eqref{neg}.  
\qquad\halmos
\vskip.1in

As a consequence of Proposition \ref{prop2} we only need consider $G(u-\RbarG_{\guG}, \RGu)$ in Theorem \ref{ThTS} when $\whL^{-1}_{T_u-}\le \tuG< \whL^{-1}_{T_u}.$  To relate this directly to the excursion $\geps$ we introduce
\be\label{Kdef}
K(z,w)=
\begin{cases}
G(z-\wbar_{g_{z}},w^{z})& \text{ if } \tz<\infty\\
1& \text{ else},
\end{cases}
\ee
and observe that if $\whL^{-1}_{T_u-}\le\guG\le \tuG< \whL^{-1}_{T_u},$ then
\be\label{RG=de}
(u-\RbarG_{\guG}, \RGu)=(u+\whH^\Gamma_{T_u-}-\gepsbar_{g^\geps_{u+\whH^\Gamma_{T_u-}}}, \geps^{u+\whH^\Gamma_{T_u-}}),
\ee
and so
\be\label{defGK}
G(u-\RbarG_{\guG}, \RGu)
=K(u+\whH^\Gamma_{T_u-},\geps).
\ee
Note also that by setting $z=u$ and $w=\geps$ in \eqref{Kdef} we have
\be
G(u-\gepsbar_{g^\geps_u}, \geps^u)=K(u,e_{T_u}),
\ee
thus
\be\label{G=K}
\lim_{u\to\infty}EG(u-\gepsbar_{g^\geps_u}, \geps^u)=\lim_{u\to\infty}EK(u,e_{T_u}).
\ee

\begin{prop}\label{cY} Assume $G$ satisfies the same conditions as in Theorem \ref{ThmGC} and $\Gamma$ is bounded away from 1.  Let
\be\label{L_G}
L_G=\lim_{u\to\infty}EG(u-\gepsbar_{g^\geps_u}, \geps^u).
\ee
Then for any $x\ge 0$
\be\ba\label{nuPuA}
\lim_{u\to\infty}\EuG[K(u+\whH^\Gamma_{T_u-},\geps);\whL^{-1}_{T_u-}\le \tau_{u-x}^\Gamma\le \tuG< \whL^{-1}_{T_u}]= L_G.
\ea\ee
\end{prop}

\noindent{\bf Proof}\ \  Fix $x\ge 0$.  Then by
\eqref{Tu1}, 
\be\ba\label{BE3}
E[K(u+\whH^\Gamma_{T_u-},\geps)&;\whL^{-1}_{T_u-}\le \tau_{u-x}^\Gamma\le \tuG< \whL^{-1}_{T_u}]\\
&=E[K(u+\whH^\Gamma_{T_u-},e_{T_u});\mcRG_{T_u-}\le u-x,\ebar_{T_u}>u+\whH^\Gamma_{T_u-}]\\
&=E\sum_tK(u+\whH^\Gamma_{t-},e_{t})I(\mcRG_{t-}\le u-x,\ebar_t>u+
\whH^\Gamma_{t-})\\
&=\int_0^\infty dt\int_z P(\mcRG_{t-}\le u-x,\whH^\Gamma_{t-}\in dz)\whn(K(u+z,w);
\tau_{u+z}<\gz)\\
&=\int_0^\infty dt\int_z P(\mcRG_{t-}\le u-x,\whH^\Gamma_{t-}\in dz)EK(u+z,e_{T_{u+z}})\,\whn(
\tau_{u+z}<\gz)
\ea\ee
where the final equality follows from \eqref{pppu}.

In \eqref{BE3} set $G\equiv 1,$  and hence  $K\equiv 1$  by \eqref{Kdef}, to obtain
\be\label{BE4}
P(\whL^{-1}_{T_u-}\le \tau_{u-x}^\Gamma\le \tuG< \whL^{-1}_{T_u})=\int_0^\infty dt\int_z P(\mcRG_{t-}\le u-x,\whH^\Gamma_{t-}\in dz)\whn(\tau_{u+z}<\gz).
\ee
Fix $L_1<L_G<L_2.$
Then for sufficiently large $u$ and all $z\ge 0$, by \eqref{G=K} and \eqref{L_G} we have
\be
L_1\le EK(u+z,e_{T_{u+z}})\le L_2.
\ee
Using these bounds in \eqref{BE3}, together with \eqref{BE4}, we have for large $u$
\be\ba\label{LL}
L_1 P(\whL^{-1}_{T_u-}\le\tau_{u-x}^\Gamma\le \tuG< \whL^{-1}_{T_u})& \le E[K(u+\whH^\Gamma_{T_u-},\geps);\whL^{-1}_{T_u-}\le\tau_{u-x}^\Gamma\le \tuG< \whL^{-1}_{T_u}]\\ &\le L_2 P(\whL^{-1}_{T_u-}\le\tau_{u-x}^\Gamma\le \tuG< \whL^{-1}_{T_u}).
\ea\ee
Divide throughout by $P(\tuG<\infty)$, let $u\to \infty$  and use \eqref{all} (this is where we use $\Gamma$ bounded away from $1$), then finally let $L_1\uparrow L_G$ and $L_2\downarrow L_G$ to get the result.
\qquad\halmos
\vskip.1in
 
 We will need the following consequence of \eqref{C};
 \be\label{Hint}
 \int_{z\ge 1} ze^{\ga z}\Pi_H(dz)<\infty,
 \ee
 where $\Pi_H$ is the L\'evy measure of $H$.
 This can be proved in a similar manner to Proposition 7.1 in \cite{G2}.  We will also need the following relationship between $\nbar$ and $\Pi_H$ given in Corollary 4.1 of \cite{G1};  for $z>0$
 \be\label{nbarH}
 \nbar(w_\gz\in dz)=\Pi_H(dz).
 \ee
 
Recall that if $\whL^{-1}_{T_u-}\le\guG\le \tuG< \whL^{-1}_{T_u}$  then \eqref{defGK} holds. The next result shows this holds asymptotically, in an appropriate sense, if just $\whL^{-1}_{T_u-}\le \tuG< \whL^{-1}_{T_u}.$

\begin{prop}\label{prop1} Assume $G$ satisfies the same conditions as in Theorem \ref{ThmGC} and $\Gamma$ is bounded away from 1.  Let
\be
Y_u=\big|G(u-\RbarG_{\guG}, \RGu)-K(u+\whH^\Gamma_{T_u-},\geps)\big|.
\ee
Then
\be\ba\label{enuPu}
\lim_{u\to\infty}\EuG(Y_u;\whL^{-1}_{T_u-}\le \tuG< \whL^{-1}_{T_u})= 0.
\ea\ee
\end{prop}

\noindent{\bf Proof}\ \  
First observe that if $\whL^{-1}_{T_u-}\le \tuG< \whL^{-1}_{T_u}$ then by \eqref{GO} and \eqref{os},
\be\label{eq1}
G(u-\RbarG_{\guG}, \RGu)\le C e^{\ga O_{u+\whH^\Gamma_{T_u-}}(\geps)}
\ee
 and by \eqref{Gbd} and \eqref{Kdef},
\be\label{eq2}
K(u+\whH^\Gamma_{T_u-},\geps)
=G(u+\whH^\Gamma_{T_u-}-\gepsbar_{g^\geps_{u+\whH^\Gamma_{T_u-}}}, \geps^{u+\whH^\Gamma_{T_u-}})
\le C e^{\ga O_{u+\whH^\Gamma_{T_u-}}(\geps)}.
\ee
For $x> 0$ set 
\be
A_{u,x}=\{u+\whH^\Gamma_{T_u-}-\gepsbar_{g^\geps_{u+\whH^\Gamma_{T_u-}}}< x\}.
\ee
If $\whL^{-1}_{T_u-}\le \tau_{u-x}^\Gamma\le \tuG< \whL^{-1}_{T_u}$ and $A_{u,x}$ hold then $\whL^{-1}_{T_u-}\le\guG\le \tuG< \whL^{-1}_{T_u},$ and so by \eqref{defGK},
$Y_u=0$.
Hence for any $x> 0$, by 
 \eqref{eq1} and \eqref{eq2}
\be\ba\label{YI,II}
\EuG(Y_u;\whL^{-1}_{T_u-}\le \tuG< \whL^{-1}_{T_u})
&\le C\EuG(e^{\ga O_{u+\whH^\Gamma_{T_u-}}(\geps)};  A^c_{u,x},\whL^{-1}_{T_u-}\le\tau_{u-x}^\Gamma\le \tuG< \whL^{-1}_{T_u})\\
&\qquad  +C\EuG(e^{\ga O_{u+\whH^\Gamma_{T_u-}}(\geps)};  \tau_{u-x}^\Gamma< \whL^{-1}_{T_u-}\le\tuG< \whL^{-1}_{T_u})\\
&=I+II.
\ea\ee

For $II,$ set
$\tG(y,w)=e^{\ga(w_\gz-y)}I(w_\gz\ge y)$ and define $\tK$ by \eqref{Kdef}.  Then on 
$\{\whL^{-1}_{T_u-}\le\tuG< \whL^{-1}_{T_u}\}$,
\be
\tK(u+\whH^\Gamma_{T_u-},\geps)=\tG(u+\whH^\Gamma_{T_u-}-\gepsbar_{g^\geps_{u+\whH^\Gamma_{T_u-}}}, \geps^{u+\whH^\Gamma_{T_u-}})=e^{\ga O_{u+\whH^\Gamma_{T_u-}}(\geps)}.
\ee
Thus by Proposition \ref{cY}, first with $x=0$ then with a general $x> 0,$ we have 
\be\ba
II
&=C\EuG(\tK(u+\whH^\Gamma_{T_u-},\geps);  \whL^{-1}_{T_u-}\le\tuG< \whL^{-1}_{T_u})\\
&\qquad -C\EuG(\tK(u+\whH^\Gamma_{T_u-},\geps);  \whL^{-1}_{T_u-}\le\tau_{u-x}^\Gamma\le \tuG< \whL^{-1}_{T_u})\\
&\to CL_{\tG}-CL_{\tG}=0.
\ea\ee

For $I$ set $\tG(y,w)=e^{\ga(w_\gz-y)}I(w_\gz\ge y)I(y\ge x)$ and again define $\tK$ by \eqref{Kdef}.  Then
by Proposition \ref{cY} 
\be\ba
I&= C\EuG[\tK(u+\whH^\Gamma_{T_u-},\geps);\whL^{-1}_{T_u-}\le\tau_{u-x}^\Gamma\le \tuG< \whL^{-1}_{T_u}]
\to CL_{\tG},
\ea\ee
where
\be\ba
L_{\tG}
&= \int_{[0,\infty) }\frac{\ga}q e^{\ga y}dy \int_{{D}}\tG(y,\gep)\ \nbar(d\gep, \gep(\gz)>y)+d_H\frac{\ga}q \tG(0,{\bf 0})\\
&=\int_{y\ge x}\frac{\ga}q e^{\ga y}dy \int_{{D}}e^{\ga(w_\gz-y)}\ \nbar(d\gep, \gep(\gz)>y)\\
&=\int_{y\ge x}\frac{\ga}q e^{\ga y}dy \int_{z>y}e^{\ga(z-y)}\Pi_H(dz)\\
&\le \int_{z\ge x}\frac{\ga}q ze^{\ga z}\Pi_H(dz)\to 0
\ea\ee
as $x\to\infty$ by \eqref{Hint}, where we used \eqref{nbarH} to obtain the final equality.
Since \eqref{YI,II} holds for any $x> 0$, we can let $u\to\infty$ then $x\to\infty$ in \eqref{YI,II}  to get \eqref{enuPu}
\qquad\halmos
\vskip.1in

\noindent{\bf Proof of Theorem \ref{ThTS}}\ \  By Propositions \ref{prop2} and \ref{prop1}
\be\ba
\big|\EuG&G(u-\RbarG_{\guG}, \RGu)-\EuG[K(u+\whH^\Gamma_{T_u-},\geps);\whL^{-1}_{T_u-}\le \tuG< \whL^{-1}_{T_u}]\big|\\
&\le \EuG[\big|G(u-\RbarG_{\guG}, \RGu)-K(u+\whH^\Gamma_{T_u-},\geps)\big|;\whL^{-1}_{T_u-}\le \tuG< \whL^{-1}_{T_u}]\\\
&\qquad +
\EuG[G(u-\RbarG_{\guG}, \RGu);\whL^{-1}_{T_u-}\le \tuG<\infty]
\to 0.
\ea\ee
The result then follows from Proposition \ref{cY} with $x=0$.
\qquad\halmos
\vskip.1in

Combining Theorem \ref{ThTS} with Proposition \ref{P31} we have

\begin{prop}\label{P41}  For any $y\ge0, x\ge 0, v\ge 0, t\ge 0$
\be\ba\label{jtconC4}
\PuG(&u-\RbarG_{\guG}\in dy,\RG_{\tuG}-u\in dx, u-\RbarG_{\guG}\in dv, \tuG- \guG\in dt)\\
&\tow\frac{\ga }q e^{\ga y} dy I(v\ge y) \whV(dt,dv-y)\Pi_X(v+dx) + d_H\frac{\ga }q \delta_{(0,0,0,0)}(dx,dy,dv,dt).
\ea\ee
\end{prop}

We can also compute the limiting expected value for certain unbounded functionals in \eqref{nuPu}.  As an example the following result gives the future value, at time $\guG$, of a Gerber-Shiu expected discounted penalty function (EDPF).

\begin{prop}\label{P42}
For any $\gl\ge 0, \eta\le \ga$ and $\gd\ge 0$, if $\eta+\gl-\ga\neq 0$
\be\ba\label{EDPFG}
\lim_{u\to\infty}\EuG e^{-\lambda(u-\RbarG_{\guG})+\eta(\RG_{\tuG}-u)-\delta(\tuG- \guG)}
&=
\lim_{u\to\infty}E e^{-\lambda(u-\gepsbar_{g_u^\geps})+\eta(\geps_{\tu^\geps}-u)-\delta(\tu^\geps- g_u^\geps)}\\&=\frac{\ga(\gk(\gd,\gl-\ga)-\gk(\gd,-\eta))}{q(\eta+\gl-\ga)}.
\ea\ee
\end{prop}
This follows from the analogous result for $X$ under $\Pu$ given in  (9.5) of \cite{G1}.  

We conclude by computing the expected present value of tax paid conditional on ruin occurring.

\begin{thm}\label{Div}  If $\int_0^\infty Ee^{-\ga \whH_t^\Gamma}dt<\infty$ 
and $\gd\ge 0$ then
\be\label{div1}
\lim_{u\to\infty} \EuG \left(\int_0^{\tuG} e^{-\gd s}\Gamma_s d|\Xubar|_s\right)= \frac{E\int_0^\infty dt\,e^{-\ga \whH_t^\Gamma}\big(\int_0^{\whL^{-1}_{t-}} e^{-\gd s}\Gamma_s d|\Xubar|_s\big)}
{\int_0^\infty Ee^{-\ga \whH_t^\Gamma}dt}.
\ee
The limit is finite if in addition  $\Gamma$ is bounded away from $0$.  Finally if $\Gamma\equiv \gamma\in (0,1)$ and $X$ is spectrally positive then
\be\label{div3}
\lim_{u\to\infty} \EuG \left(\int_0^{\tuG} e^{-\gd s}\Gamma_s d|\Xubar|_s\right)
= \frac{\gamma}{\ga(1-\gamma)+\whk(\gd,0)}.
\ee
\end{thm}

\noindent{\bf Proof}\ \  Since $d|\Xubar|_s$ does not assign mass to intervals of the form $(\whL^{-1}_{t-},\whL^{-1}_{t}),$ it follows that if $\whL^{-1}_{t-}\le\tuG< \whL^{-1}_{t}$ then
\be
\int_0^{\tuG} e^{-\gd s}\Gamma_s d|\Xubar|_s=\int_0^{\whL^{-1}_{t-}} e^{-\gd s}\Gamma_s d|\Xubar|_s.
\ee
Set 
\be
K_t=\int_0^{\whL^{-1}_{t-}} e^{-\gd s}\Gamma_s d|\Xubar|_s.
\ee
Then
\be\label{KK-}
K_{t}-K_{t-}=\lim_{r\uparrow t}(K_{t}-K_r)=\lim_{r\uparrow t}\int_{({\whL^{-1}_{r-}},{\whL^{-1}_{t-}}]} e^{-\gd s}\Gamma_s d|\Xubar|_s\le |\Xubar|_{\whL^{-1}_{t-}}-\lim_{r\uparrow t}|\Xubar|_{\whL^{-1}_{r-}}.
\ee
Thus by \eqref{Xlcts}, $P(K_{t-}=K_t \text{ for all  } t\in\mcT)=1.$
Also from \eqref{KK-}, we have $K_{t}-K_{t-}\le \whH_t-\whH_{t-}.$  
Hence $P(K_{t-}=K_t \text{ a.e. w.r.t. Lebesgue measure})=1,$ because $\whH$ is continuous except at its jump times which are at most countably.
Since $K_{t-}$ is clearly predictable, using these properties  with the compensation formula, we have
\be\ba
E&\left(\int_0^{\tuG} e^{-\gd s}\Gamma_s d|\Xubar|_s;\tuG<\infty\right)\\
&\qquad\qquad= E \sum_{t\in \mcT}\left(\int_0^{\whL^{-1}_{t-}} e^{-\gd s}\Gamma_s d|\Xubar|_s\right)I(\mcRG_{t-}\le u,\ \ebar_t>u+\whH^\Gamma_{t-})\\
&\qquad\qquad= E \sum_{t\in \mcT}K_{t-}I(\mcRG_{t-}\le u,\ \ebar_t>u+\whH^\Gamma_{t-})\\
&\qquad\qquad=E\int_0^\infty dtK_{t-}I(\mcRG_{t-}\le u)\whn(\tau_{u+\whH^\Gamma_{t-}}<\gz)\\
&\qquad\qquad=E\int_0^\infty dt\left(\int_0^{\whL^{-1}_{t-}} e^{-\gd s}\Gamma_s d|\Xubar|_s\right)I(\mcRG_{t-}\le u)\whn(\tau_{u+\whH^\Gamma_{t-}}<\gz).
\ea\ee
Dividing by $P(\tuG<\infty)$ and using \eqref{tuglim} and  \eqref{MP3}, we obtain \eqref{div1} provided we can justify taking the limit inside the integral.  If the numerator on the right side of \eqref{div1} is finite, then since 
\be
\frac{\whn(\tau_{u+\whH^\Gamma_{t-}}<\gz)}{P(\tuG<\infty)}\le Ce^{-\ga\whH^\Gamma_{t-}}
\frac{\whn(\tau_{u}<\gz)}{P(\tuG<\infty)}
\ee
for $u\ge 1$ by \eqref{npot}, we can apply dominated convergence to obtain the result.  If it is infinite then the result follows from Fatou.  

Since the numerator of the limit in \eqref{div1} is increasing in $\Gamma,$ 
to check the limit is finite when 
$\Gamma$ is bounded away from $1,$ it suffices to check the numerator is finite when 
$\Gamma\equiv \gamma\in (0,1)$.  In that case since $e^{-\gd s}\Gamma_s\le 1,$
\be\ba
E\int_0^\infty dt\,e^{-\ga \whH_t^\Gamma}\left(\int_0^{\whL^{-1}_{t-}} e^{-\gd s}\Gamma_s d|\Xubar|_s\right)
&\le  \int_0^\infty E\left(\whH_te^{-\ga(1-\gamma) \whH_t}\right)dt<\infty.
\ea\ee

In the case $X$ is spectrally positive, we may take $\whL_t=|\Xubar|_t$ in which case $\whH_t=t$.  Thus if that $\Gamma\equiv \gamma\in (0,1)$, then 
\be\label{denom}
\int_0^\infty Ee^{-\ga \whH_t^\Gamma}dt=\int_0^\infty e^{-\ga(1-\gamma)t}dt=\frac 1{\ga(1-\gamma)},
\ee
and after a change of variable
\be\ba
E\int_0^\infty dt\,e^{-\ga \whH_t^\Gamma}\left(\int_0^{\whL^{-1}_{t-}} e^{-\gd s}\Gamma_s d|\Xubar|_s\right)
&=\int_0^\infty dt\,e^{-\ga(1-\gamma)t}\int_0^{t} Ee^{-\gd \whL^{-1}_r}\gamma dr\\
&=\frac{\gamma}{\whk(\gd,0)}\int_0^\infty dt\,e^{-\ga(1-\gamma)t}(1-e^{-\whk(\gd,0)t}) \\
&=\frac{\gamma}{\ga(1-\gamma)(\ga(1-\gamma)+\whk(\gd,0))},
\ea\ee
which results in \eqref{div3} after dividing by \eqref{denom}.
\qquad\halmos

\begin{ex} We  compute the expected total tax paid ($\gd=0$) in Example \ref{eg1}.  First observe that after a change of variable
\be\ba
\int_0^{\whL^{-1}_{t-}} \Gamma_s d|\Xubar|_s&=\int_0^t f(s)ds\\
&=
\begin{cases}
0, & t\le \gb\\
t-\gb-\gb\ln(t/\gb), & t\ge  \gb.
\end{cases}
\ea\ee
Thus from \eqref{HG1}
\be\ba
E\int_0^\infty dt\,e^{-\ga \whH_t^\Gamma}\left(\int_0^{\whL^{-1}_{t-}} e^{-\gd s}\Gamma_s d|\Xubar|_s\right)
&=\int_\gb^\infty e^{-\ga\gb(1+\ln(t/\gb)}\big(t-\gb-\gb\ln(t/\gb)\big)dt\\
&=\gb^2e^{-\ga\gb}\int_1^\infty \frac{s-1-\ln s}{s^{\ga\gb}}ds\\
&=\frac{\gb^2e^{-\ga\gb}}{(\ga\gb-1)^2(\ga\gb-2)}
\ea\ee 
if $\ga\gb>2$ and is infinite otherwise.  Hence from \eqref{HG2} the limiting expected total tax paid conditional on ruin occurring, given by \eqref{div1}, is
\be
\frac{E\int_0^\infty dt\,e^{-\ga \whH_t^\Gamma}\big(\int_0^{\whL^{-1}_{t-}} \Gamma_s d|\Xubar|_s\big)}
{\int_0^\infty dt\,Ee^{-\ga \whH_t^\Gamma}}
=\frac{\ga\gb^2e^{-\ga\gb}}{(\ga\gb-1)(\ga\gb-2)(\ga\gb-1+e^{-\ga\gb})}
\ee
when $\ga\gb>2$ and is infinite otherwise.
\end{ex}



\end{document}